\documentclass{amsart}   

\usepackage{amsmath}
\usepackage{amssymb}
\usepackage{rotating}
\usepackage{graphicx}
\usepackage{url}
\usepackage{color}

\numberwithin{equation}{section}

\usepackage{amsthm}
\theoremstyle{plain}
\newtheorem{theorem}{Theorem}[section]
\newtheorem{asymptoticformula}[theorem]{Pairing-friendly curves estimate}
\newtheorem{variableDformula}[theorem]{Variable $D$ estimate}
\newtheorem{lemma}[theorem]{Lemma}

\theoremstyle{remark}
\newtheorem{remark}[theorem]{Remark}

\numberwithin{equation}{section}

\newcommand\case[1]{{\ensuremath{\mathrm(\mathit{#1}\mathrm)}}}

\newcommand\bZ{\ensuremath{\mathbb{Z}}}
\newcommand\bQ{\ensuremath{\mathbb{Q}}}
\newcommand\bR{\ensuremath{\mathbb{R}}}

\newcommand\bF{\ensuremath{\mathbb{F}}}
\newcommand\cN{\ensuremath{\mathcal{N}}}
\newcommand\ve{\ensuremath{\varepsilon}}
\newcommand\N{\ensuremath{\mathrm{N}}}
\newcommand\gp{\ensuremath{\mathfrak{p}}}
\newcommand\gr{\ensuremath{\mathfrak{r}}}
\newcommand\revision[1]{#1}

\begin{document}
\title{Heuristics on pairing-friendly elliptic curves}

\date{\today}
\author[J. Boxall]{John Boxall}
\address{Laboratoire de Math\'ematiques Nicolas Oresme, CNRS -- UMR 6139, Universit\'e de Caen Basse-Normandie, boulevard mar\'echal Juin, BP 5186, 14032 Caen cedex, France}
\email{john.boxall@unicaen.fr}

\begin{abstract}We present a heuristic asymptotic formula as $x\rightarrow \infty$ for the number of isogeny classes of pairing-friendly elliptic curves over prime fields with fixed embedding degree $k\geq 3$, with fixed discriminant, with rho-value bounded by a fixed $\rho_0$ such that $1<\rho_0<2$, and with prime \revision{subgroup} order at most $x$.\end{abstract}



\thanks{This paper was written while the author participated in the project \textit{Pairings and Advances in Cryptology for E-cash} (PACE) funded by the ANR. I would like to thank David Gruenewald for his criticisms and suggestions concerning a preliminary version of this manuscript and Igor Shparlinski for drawing my attention to \cite{ULS}.}

\maketitle

\section*{Introduction}

Pairing-based cryptography protocols first became important with the work of Joux \cite{Jx} and nowadays have numerous applications to the security of information transmission and other fields. Many of these protocols require the construction of elliptic curves over finite fields having very special properties. More precisely, let $q=p^f$ be a power of the prime $p$ and let $k\geq 1$, $r\geq 1$ be integers. We need to be able to construct an elliptic curve $E$ over the finite field $\bF_q$ with $q$ elements that satisfies the following:

\case{a} $E$ has a point $P$ of order $r$ rational over $\bF_q$;

\case{b} The group of points $E[r]$ of order $r$ of $E$ is isomorphic to $(\bZ/r\bZ)^2$ and all the points of $E[r]$ are rational over the extension field $\bF_{q^k}$ of degree $k$ of $\bF_q$.

In practical applications, if a security level of $s$ bits is required, it is generally \revision{recommended} that the integer $r$ should have at least $2s$ bits \revision{(see for example Table~1 in \cite{FST})}.  \revision{This is because the Pollard-rho algorithm is generally believed to be the best attack on the elliptic discrete logarithm problem}. The subgroup of $E(\bF_q)$ generated by $P$ should be of small index in $E(\bF_q)$. Since $\sharp(E(\bF_q))\in [(\sqrt{q}-1)^2,(\sqrt{q}+1)^2]$, so that $\sharp(E(\bF_q))\approx q$, a convenient measure of the suitability of the curve is the so-called rho-value, defined by $\rho=\frac{\log{q}}{\log{r}}$, which ideally should be close to $1$. On the other hand, the integer $k$ needs to be sufficiently small to allow efficient \revision{arithmetic} in $\bF_{q^k}$, which in practice implies that $k$ is at most about $50$. These constraints on $\rho$ and $k$ imply very strong restrictions on the choice of elliptic curve, making suitable curves very rare (\cite{BaK},  \cite{GMV}, \cite{LS}, \cite{ULS}). For this reason, a systematic search to obtain curves having parameters of cryptographic interest is completely out of the question. 

Although there is considerable recent interest in protocols where the group order $r$ is composite (\cite{BGN}, \cite{BRS}, \cite{Fr}), we \revision{shall} be concerned in this paper with the more familiar situation where $r$ is a prime number, which is assumed to be the case from now on. Since known attacks on such protocols are based on the discrete logarithm in the subgroup of order $r$ of the multiplicative group $\bF_{q^k}^\times$, \revision{and this is believed to be the same difficulty as the discrete logarithm in $\bF_{q^k}^\times$ itself}, $k$ cannot be too small. In what follows, therefore, we shall often suppose that $k\geq 3$. 

Let $E$ be an elliptic curve over $\bF_q$ satisfying \case{a}, where $r$ is a prime different from $p$. Following what has become standard usage, the smallest integer $k$ such that $q^k\equiv 1Ê\pmod{r}$ is called the embedding degree of $(E,P)$ (or just of $E$ if \revision{there} is no possibility of confusion). Alternatively, the embedding degree is just the order of $q$ in $(\bZ/r\bZ)^\times$. An argument using the characteristic polynomial of the Frobenius endomorphism (see \cite{BaK}~Theorem 1) shows that if $E$ is an elliptic curve over $\bF_q$ that satisfies \case{a} and if the embedding degree $k$ of $E$ is at least $2$, then $E$ also satisfies \case{b}.  Let $\Phi_k(w)\in \bZ[w]$ denote the $k^{\textit{th}}$ cyclotomic polynomial. Then $r$ divides $\Phi_k(q)$. On the other hand, if $t$ denotes the trace of the Frobenius endomorphism of $E$ over $\bF_q$, then $\sharp(E(\bF_q))=q+1-t$ and so $q\equiv t-1\pmod{r}$. It follows that $r$ divides $\Phi_k(q)$ if and only if $r$ divides $\Phi_k(t-1)$. Furthermore, we know \revision{from Hasse's bound} that $|t|\leq 2\sqrt{q}$ and, if we suppose in addition that $p$ does not divide $t$, then $E$ is ordinary and there exists a unique square-free positive integer $D$ and a unique integer $y>0$ such that $t^2+Dy^2=4q$. The endomorphism ring of $E$ is then an order in the imaginary quadratic field $\bQ(\sqrt{-D})$. Conversely, if $t$, $D$, $y$ are integers and if $D>0$ is square-free, $t^2+Dy^2=4q$ with $q=p^f$ a power of the prime $p$ and $p$ does not divide $t$, then a theorem of Deuring \cite{De} implies that there exists an elliptic curve $E$ over $\bF_q$ such that $\sharp(E(\bF_q))=q+1-t$. If, further, $r$ is a prime dividing both $q+1-t$ and $\Phi_k(t-1)$, and if the rho-value $\frac{\log{q}}{\log{r}}$ is close to $1$, then $E$ is suitable for pairing-based cryptography. Since we only know how to construct the curve $E$ corresponding to a choice of parameters $(t,D,y)$ when $D$ is fairly small ($D\leq 10^{15}$, say, see \cite{ES}), we shall suppose except in the last section that $D$ is fixed.  

The purpose of this note is to discuss the following \revision{heuristic asymptotic estimate}.

\begin{asymptoticformula}\label{asymptoticformula}
Let $k\geq 3$ be an integer, let $D\geq 1$ be a square-free integer and let $\rho_0\in \bR$ with $1<\rho_0<2$. \revision{We suppose that} 

\begin{enumerate}
\item\revision{$(k,D)\neq (3,3)$, $(4,1)$, $(6,3)$,}

\item\revision{If $(k,D)$ is such that there exists a complete polynomial family $(r_0,t_0,y_0)$ with generic rho-value equal to $1$ (see remark (6) below and \S~\ref{BTfamily} for detailed definitions), then $\rho_0>1+\frac{1}{\deg{r_0}}$.} 
\end{enumerate}

\revision{Let $e(k,D)=2$ or $1$ according as to whether $\sqrt{-D}$ belongs to the field generated over $\bQ$ by the $k$-th roots of unity or not, let $w_D$ be the number of roots of unity in the imaginary quadratic field $\bQ(\sqrt{-D})$ and let $h_D$ be the class number of $\bQ(\sqrt{-D})$.
Then the number of triples $(r,t,y)\in \bZ^3$ with $2\leq r\leq x$ a prime number dividing $\Phi_k(t-1)$, $t^2+Dy^2=4p$ with $p$ prime, $y>0$, $r$ dividing $p+1-t$, and $p\leq r^{\rho_0}$ is asymptotically equivalent as $x\rightarrow \infty$ to
\begin{equation}
\frac{e(k,D)w_D}{2\rho_0h_D}\int_{2}^x\frac{du}{u^{2-\rho_0}(\log{u})^2}.
\end{equation} }  
\end{asymptoticformula} 

Several remarks are in order.

\revision{(1) If $f$ is a function that is strictly positive for sufficiently large real $x$ and if $g$ is a second function defined for sufficiently large real $x$ we say that $g$ is asymptotically equivalent to $f$ as $x\rightarrow \infty$ if $g(x)=f(x)(1+o(1))$.}

\revision{(2)} Integrating by parts, we find
\begin{equation*}
\int_{2}^x\frac{du}{u^{2-\rho_0}(\log{u})^2}=\frac{1}{\rho_0-1}\frac{x^{\rho_0-1}}{(\log{x})^2}\Big(1+O\big(\frac{1}{\log{x}}\big)\Big),
\end{equation*}
where the constant implied by the $O$ is independent of $\rho_0$. Thus, for fixed $\rho_0$, the number of triples is also asymptotically equivalent to
\begin{equation}\label{ourasympt}
\frac{e(k,D)w_D}{2\rho_0(\rho_0-1)h_D}\frac{x^{\rho_0-1}}{(\log{x})^2}.
\end{equation}
However, in view of the term $\rho_0-1$ that appears in the denominator in this formula, the version with the integral seems preferable.

\revision{(3)} Several papers have appeared in the literature showing (either heuristically or unconditionally) that pairing-friendly elliptic curves are sparse (see for example \cite{BaK},  \cite{GMV} \S 4.1, \cite{LS} and \cite{ULS} and also Remark~\ref{ULSrem}). However, to the best of our knowledge, this paper is the first to suggest a possible asymptotic formula.

{(4)} One knows that two elliptic curves $E_1$ and $E_2$ over $\bF_q$ are isogenous if and only if $\sharp(E_1(\bF_q))=\sharp(E_2(\bF_q))$. It follows that to each triple there corresponds a unique isogeny class of elliptic curves, and it is clear that the embedding degree $k$ and the rho-value $\frac{\log{p}}{\log{r}}$ are invariant under isogeny. Thus (\ref{asymptoticformula}) can be interpreted as counting isogeny classes of pairing-friendly elliptic curves. For given $D$, the methods of \cite{ES} construct curves whose endomorphism ring is the maximal order of $\bQ(\sqrt{-D})$. On the other hand, Theorem~6.1 of \cite{W} shows that every isogeny class of ordinary elliptic curves contains a curve whose endomorphism ring is the maximal order of $\bQ(\sqrt{-D})$. Thus, if $D$ is sufficiently small, the methods of \cite{ES} enable one to construct at least one member of an isogeny class corresponding to any triple $(r,t,y)$.  

\revision{(5)} We have supposed that $t^2+Dy^2=4p$ with $p$ prime rather than a power of a prime. However, as is usually the case in analytic number-theoretical situations, we expect solutions with $t^2+Dy^2=4p^f$ and $f> 1$ to be negligible in number as compared with those with $f=1$, so they should not affect the asymptotic estimate. Since only finitely many primes $r$ divide $\Phi_k(-1)$, we can suppose that $t\neq 0$, in which case Deuring's theorem implies that every choice of triple $(r,t,y)$ with the properties indicated in (\ref{asymptoticformula}) corresponds to an isogeny class of ordinary elliptic curves suitable for pairing-based cryptography provided $\rho_0$ is chosen sufficiently close to $1$.   

\revision{(6)} \revision{We know of only one pair $(k,D)$ for which there is a complete polynomial family $(r_0,t_0,y_0)$ with generic rho-value equal to $1$. This is the pair $(12,3)$, and the corresponding family is the well-known Barreto-Naehrig family \cite{BaN}. In this case the degree $\deg{r_0}$ of the polynomial $r_0$ is $4$. In general, as we shall explain in \S~\ref{BTfamily}, the Bateman and Horn heuristic asymptotic formula \cite{BH} predicts that a complete polynomial family with generic rho-value equal to one will produce more triples than predicted by (\ref{asymptoticformula}) when $\rho_0<1+\frac{1}{\deg{r_0}}$. This will be a consequence of Theorem~\ref{BHtheorem} below.}

\revision{(7)} On the other hand, the cases $(k,D)=(3,3)$, $(6,3)$ and $(4,1)$ have to be excluded for a trivial reason. These are exactly the values of $(k,D)$ with $k\geq 3$ and $\bQ(\sqrt{-D})$ is equal to the field generated over $\bQ$ by the $k$-th roots of unity; one deduces easily that $t^2+Dy^2$ cannot be of the form $4p$ with $p$ a prime. See Remark~\ref{excludedkD} for further details. Recall however that this does \emph{not} imply that there are no pairing-friendly curves when $(k,D)$ takes one of these values, but only that such curves cannot be rational over prime fields. Indeed, when $(k,D)=(3,3)$, there is a well-known construction of curves over fields of square \revision{cardinality} (see \cite{FST}, \S~3.3 and also Remark~\ref{ULSrem} below).

\revision{(8)} We have excluded the cases $k=1$ and $k=2$. 

When $k=1$ and $E$ has a point $P$ of order $r$ rational over $\bF_q$, there are two possibilities:

\case{a} either all the points of $E[r]$ are rational over $\bF_q$, in which case $r^2\leq q+1+2\sqrt{q}$ by the Weil bound, which implies that the rho-value is at asymptotically least $2$, or

\case{b} the points of $E[r]$ that are not multiples of $P$ become rational only after extension of scalars to $\bF_{q^r}$, so that computations of any sort are completely infeasible. 

When $k=2$ and $E$ has a point $P$ of order $r$ rational over $\bF_q$, then $r$ divides $q+1-t$ and also $r$ divides $q+1$, since $\Phi_2(w)=w+1$. Hence $r$ divides $t$ and again there are two possibilities:

\case{a} if $t\neq 0$, then $r\leq |t|\leq 2\sqrt{q}$ and so the rho-value is asymptotically at least $2$, or

\case{b} $t=0$, in which case $E$ is supersingular. Suppose for example that the prime $r$ is such that $2r-1$ is also prime and take $q=p=2r-1$. By Deuring's theorem, there exists a supersingular elliptic curve $E$ over $\bF_p$ with $\sharp(E(\bF_p))=p+1=2r$. By the Bateman-Horn \revision{heuristics}, there is a constant $C>0$ such that number of primes $r\leq x$ with $2r-1$ prime is asymptotically equal to $C\int_{2}^x{\frac{du}{(\log{u})^2}}$. For the corresponding elliptic curves, the rho-value approaches $1$ as $r\to \infty$. Thus, when $k=2$, we expect far more pairing-friendly elliptic curves with $r\leq x$ than predicted by (\ref{asymptoticformula}).

Thus, we do not expect (\ref{asymptoticformula}) to give a reasonable estimate for the number of pairing-friendly elliptic curves when $k=1$ or $k=2$. Roughly speaking, our heuristic argument will fail in these cases because $k\in \{1,2\}$ when and only when $\Phi_k(w)$ is of degree one, and so has only the {\lq\lq}constant{\rq\rq} root $1$ or $-1\pmod{r}$ when $r$ varies. \revision{But, in view of Lemma~\ref{lemmacyc} below, it is reasonable to assume, when $k\geq 3$, that the probability that a random integer is a root of $\Phi_k(w)\mod{r}$ is $\frac{\phi(k)}{r}$.}

Here is a brief outline of the paper. In \S~\ref{heuristicarg}, we briefly describe a heuristic argument which leads to (\ref{asymptoticformula}) and in \S~\ref{numevid} we present numerical evidence for several values of $(k,D)\neq (12,3)$. In \S~\ref{BTfamily}, we review families of pairing friendly curves and in particular the Barreto-Naehrig complete family \revision{\cite{BaN}}, and explain why (\ref{asymptoticformula}) is expected to fail when \revision{$(k,D)$ satisfies condition (ii) of (\ref{asymptoticformula}) and, in particular, when  $(k,D)=(12,3)$. This involves the Bateman-Horn heuristic asymptotic estimate on polynomials with integer coefficients and its generalisation by K.~Conrad \cite{K} to polynomials with rational coefficients that take integer values.} Finally, in \S~\ref{Dvariable}, we briefly discuss a variant of (\ref{asymptoticformula}) where $D$ is allowed to vary and compare this with the recent work of Urroz, Luca and Shparlinski \cite{ULS} (see Remark~\ref{ULSrem}).  

We insist on the fact that (\ref{asymptoticformula}) is only a heuristic assertion, not a theorem. Indeed, proofs of most of the hypotheses that are used to derive it and described in \S~\ref{heuristicarg} seem to be a long way off.

All calculations reported on in this paper where done using PARI/GP \cite{pari} running on the GMP kernel \cite{GMP} and often using PARI's GP to C compiler gp2c .

\section{A heuristic argument}\label{heuristicarg}

As in the Introduction, we fix an integer $k\geq 1$ and a square-free integer $D\geq 1$. If $r$ is a prime such that $r$ does not divide $kD$, $r\equiv 1\pmod{k}$ and $-D$ is a square $\pmod{r}$, the Cocks-Pinch method \cite{CP}, as explained say in Theorem~4.1 of \cite{FST}, produces all parameters $(r,t,y)$ corresponding to ordinary curves with embedding degree $k$ and endomorphism ring an order in $\bQ(\sqrt{-D})$ having a point of order $r$. This means that $r$ divides $\Phi_k(t-1)$, $y>0$ and $t^2+Dy^2=4p$ with $p$ prime, the corresponding curve having coefficients in $\bF_p$. As is well-known, the rho-value of the curve is usually \revision{around $2$. The heuristic argument that follows will give a measure of the frequency with which it can be expected to give curves with smaller rho-values}. In what follows, we fix a real number $\rho_0$ with $1<\rho_0<2$. We wish to estimate asymptotically as $x\to \infty$ the number of triples $(r,t,y)\in \bZ^3$ as above with $r\leq x$ and $p\leq r^{\rho_0}$. Thus, the heuristic argument that follows is, in fact, an estimate of the expected number of curves with $r\leq x$ and $p\leq r^{\rho_0}$ that the Cocks-Pinch method produces. 

We first recall the following well-known Lemma, which can be extracted from \cite{Wa}, Chapter 2 \S 2:

\begin{lemma}\label{lemmacyc}
Let $k\geq 1$ be an integer and let $r$ be a prime number not dividing $k$. The following statements are equivalent.

\case{i} The cyclotomic polynomial $\Phi_k(w)$ has a root $\pmod{r}$;

\case{ii} $\Phi_k(w)$ splits into distinct linear factors $\pmod{r}$;

\case{iii} $r\equiv 1\pmod{k}$.

\case{iv} $r$ splits completely in the cyclotomic field $\bQ(\zeta_k)$ generated over $\bQ$ by a primitive $k^{\textit{th}}$ root of unity $\zeta_k$. 
\end{lemma} 

\revision{Let $r\geq 2$ be any integer. By Lemma~\ref{lemmacyc}, the probability that $r$ is prime and splits completely in $\bQ(\zeta_k)$ is equal to the probability that $r$ is prime and that $r\equiv 1\pmod{k}$. Since there are $\phi(k)$ residue classes $\pmod{k}$ consisting of integers prime to $k$, the prime number theorem generalized to arithmetic progressions implies that this is equal to $\frac{1}{\phi(k)\log{r}}$.} 

\revision{On the other hand, if $t$ is an arbitrary integer, we assume that the probability that $\Phi_k(t-1)\equiv 0\pmod{r}$ is $\frac{\phi(k)}{r}$. Since $\Phi_1(w)=w+1$ and $\Phi_2(w)=w-1$, this is reasonable only when $k\geq 3$. Thus, the probability that $r$ is prime and divides $\Phi_k(t-1)$ is $\frac{1}{\phi(k)\log{r}}\frac{\phi(k)}{r}=\frac{1}{r\log{r}}$.}

\revision{Next, we estimate the probability that $p$ be prime. To do this, we consider the element
\begin{equation*}
\pi=\frac{t+y\sqrt{-D}}{4}
\end{equation*}
of the imaginary quadratic field $\bQ(\sqrt{-D})$. Then $\pi$ is a root of $x^2-tx+p$, so that $\pi$ is an algebraic integer. Write $\N(\alpha)$ for the norm down to $\bQ$ of an element $\alpha$ of $\bQ(\sqrt{-D})$. Then $\N(\pi)=p$ so that the condition that $p$ be prime is equivalent to the condition that $\pi$ generate a principal prime ideal of $\bQ(\sqrt{-D})$. By the prime ideal theorem in $\bQ(\sqrt{-D})$ (see for example \cite{Na}, Chapter 7 \S 2), the number of principal prime ideals $\gp$ of $\bQ(\sqrt{-D})$ of prime norm $p$ bounded by $X$ is equivalent to $\frac{X}{h_D\log{X}}$ as $X\to \infty$. Applying this with $X=r^{\rho_0}$ and observing that every non-zero principal ideal of $\bQ(\sqrt{-D})$ has $w_D$ generators all having the same norm, we deduce that the expected number of primes $p\leq r^{\rho_0}$ associated to a triple $(r,t,y)$ is prime is equal to $\frac{w_Dx^{\rho_0}}{h_D\rho_0\log{r}}$.}

\revision{Finally, we estimate the probability that $r$ divides $p+1-t$, given that $r$ is prime. Now $p+1-t=\N(\pi-1)$, so that the $r$ divides $p+1-t$ if and only if there exists a prime ideal $\gr$ lying above $r$ and dividing $\pi-1$. Since $\rho_0<2$, this implies that $r$ splits in $\bQ(\sqrt{-D})$ as a product $\gr\bar{\gr}$ of two prime ideals of degree one. The probability that a random algebraic integer $\pi$ satisfies $\pi\equiv 1\pmod{\gr}$ is $\frac{1}{r}$ and the generalisation to $\bQ(\sqrt{-D})$ of Dirichlet's theorem on primes in arithmetic progressions implies that this remains true if $\pi$ generates a prime ideal. Since there are two primes ideals $\gr$ and $\bar{\gr}$ dividing $r$, the probability that $r$ divides $p+1-t$ given that it splits in $\bQ(\sqrt{-D})$ is $\frac{2}{r}$.}

\revision{On the other hand, the probability that $r$ splits as a product of two degree one primes in $\bQ(\sqrt{-D})$ is $1$ if $\sqrt{-D}\in \bQ(\zeta_k)$, and $\frac{1}{2}$ if not. This is equal to $\frac{e(k,D)}{2}$.} 

\revision{Taking all this into account and making various obvious independence hypotheses, we obtain that the number of triples $(r,t,y)$ such that $r\leq x$ is prime, $r\equiv 1\pmod{k}$, $r$ divides $\Phi_k(t-1)$, and $t^2+Dy^2=4p$ with $p\leq r^{\rho_0}$ a prime should be equivalent to
\begin{equation*}
\sum_{2\leq r\leq x}{\frac{1}{r\log{r}}\frac{w_Dx^{\rho_0}}{h_D\rho_0\log{r}}\frac{2}{r}\frac{e(k,D)}{2} } =\frac{e(k,D)w_D}{\rho_0h_D}\sum_{2\leq r\leq x}\frac{1}{r^{2-\rho_0}(\log{r})^2}.
\end{equation*} 
Here the sums are over all integers $r$ such that $2\leq r\leq x$. Since 
\begin{equation*}
\sum_{2\leq r\leq x}\frac{1}{r^{2-\rho_0}(\log{r})^2}\sim \int_{2}^x\frac{du}{u^{2-\rho_0}(\log{u})^2},
\end{equation*}   
this estimate differs by a factor of $2$ from that in (\ref{asymptoticformula}), the difference being due to the fact that we assumed in (\ref{asymptoticformula}) that $y>0$ whereas in the preceding argument the sign of $y$ is arbitrary}.

\begin{remark} \label{excludedkD}
The independence hypotheses alluded to above \revision{assume} that  \revision{$\pi$} is an essentially random element of the set of algebraic integers of $\bQ(\sqrt{-D})$ such that \revision{$\pi-1$} belongs to one of the prime ideals dividing $r$. In particular, the probability that it generates a prime ideal should be that predicted by the prime ideal theorem. This is not true when $(k,D)=(3,3)$, $(6,3)$ or $(4,1)$, in other words in those cases where $\bQ(\zeta_k)=\bQ(\sqrt{-D})$. Suppose for example that $(k,D)=(3,3)$. The condition $r|\Phi_3(t-1)$ then implies that $4r$ divides $4t^2-4t+4$. On the other hand, since $4r$ divides $(t-2)^2+3y^2=t^2-4t+4+3y^2$, we find by subtraction that $4r$ divides $3(t^2-y^2)$. When $r\geq 5$, this implies that $t\equiv \pm y\pmod{4r}$. Since $|t|\leq 2r$ and $|y|\leq 2r$, this implies that $t=\pm y$ when $r$ is sufficiently large and so $t^2+3y^2$ cannot be of the of the form $4p$ with $p$ a prime. A similar argument works when $(k,D)=(6,3)$ or $(4,1)$. Thus the use of the prime ideal theorem is not justified in these cases. 
\end{remark}

\section{Numerical evidence}\label{numevid}

In order to test (\ref{asymptoticformula}) numerically, we wrote a programme in PARI/GP \cite{pari} to search for all triples $(r,t,y)$ with $r$ in some interval $[a,b]$, $k$, $D$ and $\rho_0$ being given. Thus for each prime $r\equiv 1\pmod{k}$ belonging to $[a,b]$ such that $-D$ is a square $\pmod{r}$, the programme finds all the roots of $\Phi_k(t-1)\equiv 0\pmod{r}$, searches for those for which \revision{$|t|\leq 2r^{\frac{\rho_0}{2}}$} and then those for which there exists $y>0$ such that $t^2+Dy^2=4p$ with $p$ prime and $p\leq r^{\rho_0}$, and outputs the vector of all sextuples $(r,t,y,h,p,\rho)$ with $r$, $t$, $y$ and $p$ as before, \revision{$h$ the} cofactor defined by $p+1-t=rh$, and \revision{$\rho=\frac{\log{p}}{\log{r}}$ the} actual rho-value. 

For a given $r$, there are two possible strategies for finding $t$. The first is to factor $\Phi_k(x)\pmod{r}$ using a standard factorisation algorithm for univariate polynomials over finite fields. The second is to first choose at random a primitive root $g \pmod{r}$, so that if $s=g^{\frac{r-1}{k}}\pmod{r}$, then $s$ is a primitive $k$-th root of unity in the field with $r$ elements. The possible values of $t$ are then $s^\ell+1\pmod{r}$ as $\ell$ ranges over the integers between $1$ and $k$ that are prime to $k$. This is justified by the fact the roots of $\Phi_k$ are precisely the primitive $k$-th roots of unity. In the range where the systematic search for all triples $(r,t,y)$ is feasible, the second method turned out to be the faster although it is clear that for large values of $r$ the first method is preferable since $k\leq 50$ and the exponentiation to the power $\frac{r-1}{k}$ becomes costly.

In view of the discussion in \S~\ref{heuristicarg}, our programme is basically an implementation of the Cocks-Pinch method that selects only those curves with $\rho\leq \rho_0$. However, as all primes $r\equiv 1\pmod{k}$ need to be tested, this cannot be expected in reasonable time to find curves in an interval $[a,b]$ where $a$ and $b$ are of a sufficiently large size for the curves to be of cryptographic interest (unless the value $\rho$ is taken to be close to $2$). \revision{In practice, it was found that for given $k$ and $D$ the vector of all sextuples $(r,t,y,h,p,\rho)$ could be calculated in between $15$ and $75$ seconds when $b-a=10^8$ and $b$ is smaller than about $10^{15}$. Under these conditions, the} time taken was roughly proportional to $1/\phi(k)$.  Also, in view of the irregularity that one expects when $k$ and $D$ vary and $r$ is very small, it was decided to restrict attention to $r\geq 10^6$.  

In what follows we present, for different values of $k$, $D$, $\rho_0$, $a$ and $b$, the number $N=N(k,D,\rho_0,a,b)$ of triples $(r,t,y)$ as in (\ref{asymptoticformula}) with $a\leq r\leq b$ and, for comparison, the value of the corresponding integral
\begin{equation}\label{defIrho}
I=I(k,D,\rho_0,a,b)=\frac{e(k,D)w_D}{2\rho_0 h_D}\int_a^b\frac{du}{u^{2-\rho_0}(\log{u})^2}.
\end{equation}
We define $I_0$ by $I_0(k,D,\rho_0,a,b)=e(k,D)^{-1}I(k,D,\rho_0,a,b)$: note that $I_0$ depends only on $D$ and $\rho_0$ but not on $k$.

For convenience, the tables of numerical data have been placed near the end of the paper.
 
Figure~\ref{tablepreone} gives the values of $N(k,D,1.7,10^6,85698768)$ for all $k$ such that $3\leq k\leq 30$ and all squarefree $D$ with $D\leq 15$ as well as $D=19$, $23$, $43$ and $47$. This choice of $D$ includes all imaginary quadratic fields of class number one except $\bQ(\sqrt{-163})$ and, for each integer $h$ less than or equal to $5$ at least one field whose class number is equal to $h$. The second line of the table recalls the class number $h_D$ of $\bQ(\sqrt{-D})$. The third line gives the value of $I_0=e(k,D)^{-1}I(k,D,1.7,10^6,85698768)$. The values of $I_0$ are \revision{the reason} for the choice of $85698768$ as upper limit. In fact, when $D$ is such that $w_D=2$ and $h_D=1$, then $I_0=\frac{1}{1.7}\int_{10^6}^{85698768}\frac{du}{u^{0.3}(\log{u})^2}\simeq 1000.00$ so that the predicted value of $N(k,D,1.7,10^6,85698768)$ is $1000$ in these cases. The main part of the table contains the values of $N(k,D,1.7,10^6,85698768)$, the entries corresponding to values of $(k,D)$ with $e(k,D)=2$ are marked with an asterisk;  (\ref{asymptoticformula}) predicts that they should be close to $2I_0$ and therefore roughly twice as large as the other entries in the same column. The last line of Figure~\ref{tablepreone} gives the average value of each column as $k$ varies from $3$ to $30$, the cases where $e(k,D)=2$ being counted with weight $\frac{1}{2}$ and the excluded values $(k,D)=(3,3)$, $(4,1)$ and $(6,3)$ omitted. (\ref{asymptoticformula}) predicts that each of these averages be close to the corresponding value of $I_0$.

Figure~\ref{tableone} gives the values of $N(k,D,1.5,10^6,2\times 10^8)$ for the same values of $(k,D)$ as Figure~\ref{tablepreone}. When $D$ is such that $w_D=2$ and $h_D=1$, \revision{we now have $I_0=\frac{1}{1.5}\int_{10^6}^{2\cdot 10^8}\frac{du}{u^{0.5}(\log{u})^2}\simeq 58.17$.}

Although all the entries in Figures~\ref{tablepreone} and \ref{tableone} (with the exception of those for $(k,D)=(3,3)$, $(4,1)$ and $(6,3)$) are of the order of magnitude predicted by (\ref{asymptoticformula}), there is considerable variation in the actual values, especially in Figure~\ref{tableone}. This is perhaps not unexpected, as similar variation occurs when one computes the number of values for which polynomials simultaneously take prime values and compares the result to the Bateman-Horn \revision{heuristics}. In fact, if $\pi(x)$ denotes as usual the number of primes less than or equal to the real positive $x$, no explicit formula analogous to Riemann's formula for $\pi(x)-\int_{2}^x{\frac{du}{\log{u}}}$ seems to be known in the Bateman-Horn context (see for example \cite{KorRie} for a discussion of the case of prime pairs). So, presumably it would also be a hard problem to find one in the context of (\ref{asymptoticformula}).  

In order to obtain numerical data for larger values of $x$ and examine what happens when $\rho_0$ varies, it is necessary to restrict the values of $k$ and $D$. The case $(k,D)=(12,3)$ will be discussed in the next section. Figure~\ref{table3kD} presents data for the three cases \revision{$(k,D)=(28,1)$}, $(27,11)$ and $(8,23)$. In each case, they give the values of $N(\rho_0)=N(k,D,\rho_0,a,b)$ and $I(\rho_0)=I(k,D,\rho_0, a,b)$ for $\rho_0\in \{1.1,1.2,1.3,1.4,1.5\}$ and for each of the three intervals $(a,b)=(10^6,10^8)$, $(10^8,10^{10})$ and $(10^{12}-10^{10},10^{12}+10^{10})$. These results emphasize just how rare triples with rho-values close to one are. For example, if one wanted to construct a table like Figure~\ref{tablepreone} with $I_0=1000$ but taking $\rho_0=1.2$ instead of $1.7$, (\ref{asymptoticformula}) suggests that one would need to test all $r$ up to about $7.9\times 10^{29}$, which is obviously completely out of the question.

\section{The Barreto-Naehrig family and the case $k=12$, $D=3$}\label{BTfamily}

The various known methods of constructing pairing-friendly elliptic curves are reviewed in \cite{FST}. Since (\ref{asymptoticformula}) is primarily concerned with ordinary elliptic curves over prime fields and assumes that $k\geq 3$, we limit our attention to those methods which apply in these situations. We want to understand asymptotically as $x\to \infty$ the number of triples $(r,t,y)$ with $r\leq x$ that belong to such families and have rho-value at most $\rho_0$ and compare this with    
the estimate in (\ref{asymptoticformula}). Clearly we can only compare constructions where $k$ and $D$ are fixed.

Apart from the Cocks-Pinch method, which constructs all parameters corresponding to ordinary curves and on which our heuristic estimate is based, the other well-known constructions with $k$ and $D$ fixed are the polynomial families. These fall into two kinds: \case{a} sparse families, of which the most familiar example is MNT families \cite{MNT}; \case{b} complete families, of which the general construction is due to Brezing and Weng \cite{BW}. We refer to \cite{FST}, \S~5 and 6 for a detailed review of the two kinds of families.

The idea behind both constructions is to find polynomials $r_0(w)$, $t_0(w)$ and $p_0(w)\in \bQ[w]$ such that $r_0(w)$ divides both $\Phi_k(t_0(w)-1)$ and $p_0(w)+1-t_0(w)$. One then seeks values $w_0$ of $w$ for which $r_0(w_0)$, $t_0(w_0)$ and $p_0(w_0)$ are all integers with $r_0(w_0)$ prime (or a prime multiplied by a very small factor) and $p_0(w_0)$ is prime (or a prime power). The values of the integral parameters $r$, $t$ and $p$ are then respectively $r_0(w_0)$, $t_0(w_0)$ and $p_0(w_0)$ with $r_0(w_0)$ and $p_0(w_0)$ prime. \revision{By definition, the generic rho value of the family is $\frac{\deg{p_0}}{\deg{r_0}}$. As $w_0$ tends to infinity, the rho-value of the elliptic curve corresponding $w_0$ approaches the generic rho-value.}

\revision{However, the} two constructions differ in the way they treat the parameter $y$. Define the polynomial $h_0(w)$ by $p_0(w)+1-t_0(w)=r_0(w)h_0(w)$. If $r=r_0(w_0)$, $t=t_0(w_0)$, $p=p_0(w_0)$ and $h=h_0(w_0)$, then the corresponding $y$ parameter satisfies
\begin{equation*}
Dy^2=4p-t^2=4hr-(t-2)^2.
\end{equation*}

In the case of sparse families, the general idea is choose $r_0$, $t_0$ and $p_0$ in such a way that $4p_0(w)-t_0(w)^2$ is of degree two. When this is the case, the affine curve with \revision{$(w,y)$}-equation $Dy^2=4p_0(w)-t_0(w)^2$ is of genus $0$. If this curve is to have infinitely many integral points, its real locus must be either a parabola or a hyperbola. In all the cases of which we are aware, the real locus is a hyperbola. Thus, an affine change of coordinates transforms this into a generalised Pell equation $Z^2-aY^2=b$, with $a>0$ is not a square. The integral solutions of this are of the form $Z+\sqrt{a}Y=\alpha\ve^n$, where $\alpha$ runs through a finite set of elements of the real quadratic field $\bQ(\sqrt{a})$, $\ve$ is a fundamental unit of $\bQ(\sqrt{a})$, and $n\in \bZ$. From this we deduce that the number of values of $r\leq x$ that can arise from a sparse family is $O((\log{x})^2)$. On the other hand, (\ref{asymptoticformula}) predicts that there \revision{are} at least $>\!\!\!> \frac{x^{\rho_0-1}}{(\log{x})^2}$ choices of the parameters $(r,t,y,p)$ with $r\leq x$ and $p\leq r^{\rho_0}$. Thus, sparse families can only contribute a negligible proportion of pairing friendly-curves with given $k$ and $D$. 

In the case of complete families, the basic strategy was described in full generality by Brezing and Weng \cite{BW}. In addition to $r_0$, $t_0$, $h_0$ and $p_0$, we also require a polynomial $y_0$ such that $t_0(w)^2+Dy_0(w)^2=4p_0(w)$, so that the $y$ parameter is the corresponding value $y_0(w_0)$. Now, the polynomials $r_0$, $t_0$, $y_0$, $h_0$, $p_0$ simultaneously take integral values at integers $w_0$ varying over a finite set of congruence classes modulo some fixed integer.  \revision{Furthermore, if $r_0$ and $p_0$ are to give rise to triples $(r,t,y)$ corresponding to elliptic curves, they must simultaneously take prime values.}

\revision{Before going further, we recall the Bateman-Horn \revision{heuristics} \cite{BH} in the case of two polynomials $f$ and $g$ with integral coefficients. We assume that $f$ and $g$ are distinct and irreducible. For any prime $p$ let $N_p$ denote the number of solutions of the congruence $f(x)g(x)\equiv 0\pmod{p}$ and suppose that $N_p<p$ for all $p$. Then let $C$ be given by the conditionally convergent infinite product
\begin{equation}\label{BHconstant}
C=\prod_{p\geq 2 \text{ prime}}\Big(1-\frac{N_p}{p}\Big)\Big(1-\frac{1}{p}\Big)^{-2}.
\end{equation}
Then the number of integers $w_0$ with $2\leq w_0\leq X$ such that $f(w_0)$ and $g(w_0)$ are simultaneously prime is asymptotically equivalent to
\begin{equation}\label{BHasymptotic}
\frac{C}{\deg{r_0}\deg{p_0}}\int_2^X\frac{du}{(\log{u})^2}
\end{equation}
as $X\to \infty$. In particular, since $C>0$, there are infinitely many $w_0$ such that $f(w_0)$ and $g(w_0)$ are simultaneously prime.}

\revision{We need to adapt this statement to polynomials whose coefficients are rational. Let $f$, $g\in \bQ[w]$ and let $n\geq 1$ be a common denominator of the coefficients of $f$ and $g$. Then there are integers $m_i$ with $0\leq m_i<n$ such that $f(nw_0+m_i)\in \bZ$ and $g(nw_0+m_i)\in \bZ$ for all $i$ and for all $w_0\in \bZ$.}

\revision{Then, for each $i$, we can apply the generalization by K.~Conrad (see \S~2 of \cite{K}) of the Bateman-Horn heuristics to the pair of polynomials $w\mapsto f(nw_0+m_i)$ and $w\mapsto f(nw_0+m_i)$. This implies that (\ref{BHasymptotic}) still holds, although the value of $C$ will no longer be given by (\ref{BHconstant}) in general, but can be computed using Conjecture~5 of \cite{K}. Since in what follows we only need the actual value of $C$ in the case of polynomials with integer coefficients, we do not discuss this in detail.} 

\revision{Returning to our discussion of complete families, it follows that there exists a constant $C'>0$ such that the number 
of triples $(r,t,y)$ with $r\leq x$ coming from the family is asymptotically equivalent to
\begin{equation}  \label{BWasympt}
\frac{C'}{\deg{r_0}\deg{p_0}}\int_{2}^{(x/c_{r_0})^{1/\deg{r_0}}}\frac{du}{(\log{u})^2}\sim \frac{C'}{c_{r_0}^{1/\deg{r_0}}\deg{r_0}\deg{p_0}}\frac{x^{1/\deg{r_0}}}{(\log{x})^2}, 
\end{equation}
where} $c_{r_0}$ is the leading coefficient of $r_0$ and $\deg{r_0}$ is the degree of $r_0$, and the asymptotic equivalence of the two displayed formulae is seen by integrating by parts. (\revision{Note that in general $C'$ will not be equal to $C$, since both positive and negative values of $w_0$ may yield triples $(r,t,y)$.}) 

As $x_0\to \infty$, the rho-value of the triple $(r_0(w_0),t_0(w_0),y_0(w_0))$ approaches $\frac{\deg{p_0}}{\deg{r_0}}$. Comparing (\ref{ourasympt}) and (\ref{BWasympt}), we deduce that
if $\frac{1}{\deg{r_0}}>\rho_0-1$, then the Bateman-Horn \revision{heuristics} implies the complete family parametrised by $r_0$, $t_0$, \dots,  asymptotically contains more choices of triples than predicted by (\ref{asymptoticformula}). \revision{On the other hand, the rho-value of the triples $(r_0(w_0),t_0(w_0),y_0(w_0))$ tends to the generic rho-value $\frac{\deg{p_0}}{\deg{r_0}}$ as $w_0\to \infty$, so that this family can contain infinitely many triples with rho-value $\leq \rho_0$ only if $\frac{\deg{p_0}}{\deg{r_0}}\leq \rho_0$. It is clear that $\deg{p_0}\geq \deg{r_0}$ so, since $\deg{p_0}$ and $\deg{r_0}$ are integers, the conditions $\frac{\deg{p_0}}{\deg{r_0}}\leq \rho_0$ and $\frac{1}{\deg{r_0}}>\rho_0-1$ are satisfied only if $\deg{p_0}=\deg{r_0}$.  We deduce (i) of the following}
\begin{theorem}\label{BHtheorem}
\revision{We keep the notation that has just been introduced and assume the Bateman-Horn heuristics together with their generalization by K.~Conrad.} 
\begin{enumerate}
\item \revision{Suppose that $\rho_0<1+\frac{1}{\deg{r_0}}$. Then the complete family $(r_0,t_0,y_0)$ asymptotically contains more choices of parameters than predicted by (\ref{asymptoticformula}). Furthermore, one has $\deg{p_0}=\deg{r_0}$.}

\item \revision{On the other hand, if $\rho_0>1+\frac{1}{\deg{r_0}}$ then the family does not contain sufficiently many triples to contradict (\ref{asymptoticformula}).}
\end{enumerate}   
\end{theorem} 
\revision{Point (ii) is proved in a similar way to (i), again comparing of (\ref{ourasympt}) and (\ref{BWasympt}).}

\revision{On the other hand, what happens when $\rho_0=1+\frac{1}{\deg{r_0}}$ depends on the relative values of the constants appearing in (\ref{ourasympt}) and the right hand side of (\ref{BWasympt}).}

Table~8.2 of \cite{FST} summarizes, for all $k$ up to $50$, the construction of the family with the smallest rho-value and the corresponding value of $D$. When $k\geq 4$, the families listed are all complete families, and all have $\deg{p_0}>\deg{r_0}$ \emph{except when $k=12$}, in which case the corresponding value of $D$ is $3$. When $k=3$, the family is also a complete family and $D=3$ and also satisfies $\deg{p_0}=\deg{r_0}$, except that $p_0(w)=(3w-1)^2$ cannot represent primes (see \S~3.3 of \cite{FST}).  

The case $k=12$ and $D=3$ is thus expected to provide a genuine counterexample to (\ref{asymptoticformula}). The corresponding family is the well-known Barreto-Naehrig family \cite{BaN},  where
\begin{align*}
r_0(w)&=36w^4+36w^3+18w^2+6w+1,\quad t_0(w)=6w^2+1, \quad  h_0(w)=1,\\ 
                         &y_0(w)=6w^2+4w+1,\qquad  p_0(w)=36w^4+36w^3+24w^2+6w+1. 
\end{align*}
Since the degree of $r_0$ is $4$, we expect the family to provide more curves than (\ref{asymptoticformula}) when $\rho_0<1.25$.

This can be tested numerically using similar calculations to those presented in \S~\ref{heuristicarg}. To see the contribution of the Barreto-Naehrig family, we need to calculate the constant $C$ appearing in the Bateman-Horn \revision{heuristics} for it. For any prime $p$, let $N_{r_0,p}$ denote the number of solutions of $r_0(w)\equiv 0\pmod{p}$ and define $N_{p_0,p}$ similarly. Write $N_p$ for the number of solutions of $r_0(w)p_0(w)\equiv 0\pmod{p}$. Then $N_2=N_3=0$ and $N_p=N_{r_0,p}+N_{p_0,p}$ when $p\geq 5$ since $p_0(w)=r_0(w)+6w^2$ so that $r_0$ and $p_0$ cannot have a common root $\pmod{p}$. Since $r_0$ and $p_0$ have integral coefficients, the Bateman-Horn constant \revision{is given by (\ref{BHconstant}).}

As written, the \revision{product (\ref{BHconstant})} is conditionally convergent and therefore unsuitable for numerical computation. Instead, we apply the formula given by the theorem of Davenport and Schinzel  \cite{DS}. This gives
\begin{equation*}
C=\frac{\gamma}{\rho(K_{r_0})\rho(K_{p_0})}\prod_{p\geq 5}\Big(1-\frac{N_p}{p}\Big)\Big(1-\frac{1}{p}\Big)^{-N_p}\prod_{p\geq 5}\Big(1-\frac{1}{p^2}\Big)^{-N_p^{(2)}}\Big(1-\frac{1}{p^4}\Big)^{-N_p^{(4)}},
\end{equation*}
where $N_p^{(2)}$ and  $N_p^{(4)}$ denote respectively the number of irreducible factors of $r_0(x)p_0(x)\pmod{p}$ of degree $2$ and of degree $4$, $\rho(K_{r_0})$ and $\rho(K_{p_0})$ the residue at $1$  of the zeta function of the number fields $K_{r_0}$ and $K_{p_0}$ generated over $\bQ$ by a root of $r_0$ and a root of $p_0$ and 
\begin{equation*}
\gamma=\Big(1-\frac{1}{2^2}\Big)^{-2}\Big(1-\frac{1}{3^2}\Big)^{-1}\Big(1-\frac{1}{3}\Big)^{-1}=3.
\end{equation*}

The two infinite products in the Davenport-Schinzel formula for $C$ are now absolutely convergent. When $p\geq 5$ the table that follows gives the value of $N_p^{(j)}$ when $j=2$ and $j=4$:
\medskip

\begin{center}
\begin{tabular}{|c|c||c|c|c|}
\hline
$p\mod{12}$&$p_0(w) \mod{p}$&$N_p$&$N_p^{(2)}$&$N_p^{(4)}$\\
\hline
$1$&$4$ roots&$8$&$0$&$0$\\
\hline
$1$&$0$ roots&$4$&$2$&$0$\\
\hline
$5$&{}&$0$&$2$&$1$\\
\hline
$7$&{}&$2$&$3$&$0$\\
\hline
$11$&{}&$0$&$4$&$0$\\
\hline
\end{tabular}
\end{center}
\medskip

Using these formulae and taking the product over all $p$ with $5\leq p\leq 10^6$, we find that the first product appearing in the formula for $C$ equals $0.88576\dots $ and the second equals $1.26250\dots $. On the other hand, $\rho(K_{r_0})=0.36105\dots$ and $\rho(K_{p_0})=0.52642\dots$. \revision{It follows that $C\simeq 17.651$. On the other hand, since neither of the polynomials $r_0$ and $p_0$ are even functions, the values of $r_0(w_0)$ and $p_0(w_0)$ at negative integers $w_0$ will, with finitely many exceptions, be different to those at positive integers.  Hence $C'=2C$ so that $\frac{C'}{16}\simeq 2.206$} and, if the Bateman-Horn \revision{heuristics} are correct, we can expect the number of triples $(r,t,y)$ arising from the Barreto-Naehrig family with $x'\leq r\leq x$ should be approximately equal to
\begin{equation*}
\revision{J_{BN}(x',x)={}}  2.206\int_{x^{\prime 1/4}/\sqrt{6}}^{x^{1/4}/\sqrt{6}}\frac{du}{(\log{u})^2}.
\end{equation*}  

The following table gives the values of $N(12,3,\rho_0,10^6,10^8)$ together with $N(12,3,\rho_0,10^8,10^{10})$ for $\rho_0\in \{1.1,1.2,1.3,1.4,1.5\}$ and compares them with the corresponding expected value of $I(12,3,\rho_0, a,b)$.
\medskip

\revision{\begin{center}
\begin{tabular}{|c||c|c|c|c|c|}
\hline
$\rho_0$      &  $1.1$&   $1.2$&   $1.3$&   $1.4$&   $1.5$\\
\hline
\hline
$N(10^6,10^8)$&   $3$&     $8$&     $21$&    $57$&    $305$\\
\hline
$I(10^6,10^8)$&  $0.49$&  $2.25$&  $10.66$& $51.58$& $255.11$\\
\hline
$N(10^8,10^{10})$&   $6$&    $10$&     $44$&      $221$&        $1655$\\
\hline
$I(10^8,10^{10})$&$0.47$& $3.43$&  $25.83$& $199.07$& $1567.0$\\
\hline 
\end{tabular}
\end{center}}
\medskip

The column $\rho_0=1.1$ of the table contains $3$ triples with $10^6\leq r\leq 10^8$ and $6$ with $10^8\leq r\leq 10^{10}$. All these nine triples $(r,t,y)$ are in fact members of the Barreto-Naehrig family: they correspond to the values of the polynomials $r_0(x)$ etc. at $x=-107$, $-55$, $-52$, $-41$, $-15$, $20$, $78$, $82$, $123$. \revision{This should be compared with the expected contributions from the Barreto-Naehrig family which are respectively $J_{BN}(10^6,10^8)=6.05$ and $J_{BN}(10^8,10^{10})=10.26$.}

\section{What happens when $D$ varies}\label{Dvariable} 

Let again $D$ denote a square-free positive integer. As before, we denote the discriminant of the imaginary quadratic field $\bQ(\sqrt{-D})$ by $d_D$, thus $d_D=-D$ if $D\equiv 3\pmod{4}$ and $d_D=-4D$ if $D\equiv 1$, $2\pmod{4}$. If $z$ is small with respect to $x$, (\ref{asymptoticformula}) suggests that the number of triples $(r,t,y)$ as above with $r\leq x$, $p\leq r^{\rho_0}$ and $|d_D|\leq z$ should be equivalent to
\begin{equation*}
\sum_{|d_D|\leq z}\frac{e(k,D)w_D}{2\rho_0 h_D}\int_{2}^{x}\frac{du}{u^{2-\rho_0}(\log{u})^2}\end{equation*}
Here we shall not try to give a precise meaning to the condition that $z$ be small with respect to $x$, which would require a discussion of the error term in (\ref{asymptoticformula}) which would take us too far afield. We content ourselves with a heuristic asymptotic estimate for the sum
\begin{equation*}
\sum_{|d_D|\leq z}\frac{e(k,D)w_D}{2\rho_0 h_D}
\end{equation*}
as $z\to \infty$. It is well-known that $\sqrt{-D}\in \bQ(\zeta_k)$ if and only if $d_D$ divides $k$. Furthermore, $w_D=2$ except when $D=1$ or $D=3$. Therefore
\begin{equation*}
\sum_{|d_D|\leq z}\frac{e(k,D)w_D}{2\rho_0 h_D}=\frac{1}{\rho_0}\sum_{|d_D|\leq z}\frac{1}{h_D}+O(1),
\end{equation*}
where the constant implied by the $O(1)$ depends only on $k$. Estimates for the sum $\sum_{|d_D|\leq z}h_D^\alpha$ for various positive values of $\alpha$, and in particular $\alpha=1$, have been studied since the time of Gauss (see for example \cite{GH} and the references cited therein). However, we have been unable to find any reference to the case $\alpha=-1$ which is of interest here. On the other hand, heuristic considerations involving the prime ideal theorem and the residue of zeta functions at $s=1$ for imaginary quadratic fields suggest that
\begin{equation*}
\sum_{|d_D|\leq z}\frac{1}{h_D}\sim \frac{6}{\pi}\sqrt{z},   \qquad z\to \infty
\end{equation*}
and this seems to be confirmed by numerical calculation. This suggest the following heuristic
\begin{variableDformula}
Let $k\geq 3$ and $\rho_0$ such that $1<\rho_0<2$ be fixed. If $z$ is small with respect to $x$, then, as $x\to \infty$ the number \revision{$\cN(k,z,\rho_0,x)$} of triples $(r,t,y)$ as in (\ref{asymptoticformula}) with $|d_D|\leq z$ is equivalent to
\begin{equation*}
\frac{6}{\rho_0\pi}\sqrt{z}\int_{2}^x{\frac{du}{u^{2-\rho_0}(\log{u})^2}}.
\end{equation*}
\end{variableDformula}
In particular, if we can take $z=x^\alpha$ for some small positive $\alpha$ then, integrating by parts, we find that the number of triples $(r,t,y)$ with $r\leq x$ and $|d_D|\leq x^\alpha$ should be equivalent to
\begin{equation*}
\frac{6}{\rho_0(\rho_0-1)\pi}\frac{x^{\frac{\alpha}{2}+\rho_0-1}}{(\log{x})^2}.
\end{equation*}
At present it is not quite clear how large we can take $\alpha$ for this estimate to be reasonable. This depends in particular on the size of the error term in (\ref{asymptoticformula}), a problem which certainly deserves study but we prefer to leave this for future work. One reason for this is that, to the best of our knowledge, no detailed discussion of the error term in the Bateman-Horn \revision{heuristics} has appeared in the literature up till now.

\begin{remark}\label{ULSrem}
In \cite{ULS}, Urroz, Luca and Shparlinski prove a result which implies an unconditional upper bound on \revision{$\cN(k,z,\rho_0,x)$}. In fact, their Theorem 1 implies that
\begin{equation*}
\revision{\cN(k,z,\rho_0,x)}<\!\!\!< \phi(k)\big(x^{\rho_0-1}+x^{\frac{\rho_0}{2}}\big)z^{\frac{1}{2}}\frac{\log{x}}{\log{\log{x}}}<\!\!\!<  \phi(k)x^{\frac{\rho_0}{2}}z^{\frac{1}{2}}\frac{\log{x}}{\log{\log{x}}},
\end{equation*}
where the constants implied by the $<\!\!\!< $ are absolute. This follows from the hypothesis that $1<\rho_0<2$, the variable $x$ of \cite{ULS} corresponds to our $x^{\rho_0}$, the $y$ of \cite{ULS} to our $x$, and the $z$ of \cite{ULS} is contained between $\frac{1}{4}z$ and $z$ when $z$ is used in our sense. For constant $z$, this is much weaker than (\ref{asymptoticformula}), but when $(k,D)=(3,3)$ there exists the complete family 
\begin{align*}
r_0(w)=9w^2-&3w+1, \qquad t_0(w)=-3w+1, \qquad y_0(w)=3w-1\\
&h_0(w)=1,\qquad q_0(w)= (3w-1)^2,
\end{align*}
together with a similar family with $r_0(w)=9w^2-9w+3$ (see \cite{FST}, \S~3.3). The Bateman-Horn \revision{heuristics} therefore implies that $\revision{\cN(3,z,\rho_0,x)}>\!\!\!> \frac{x^{\frac{1}{2}}}{(\log{x})^2}$ for any $z\geq 3$ and any $\rho_0$. A similar argument using the Barreto-Naehrig family suggests that $\revision{\cN(12,z,\rho_0,x)} >\!\!\!> \frac{x^{\frac{1}{4}}}{(\log{x})^2}$ for any $z\geq 3$ and any $\rho_0$. Thus, the Urroz-Luca-Shparlinski upper bound for a given $k$ is strongly related to the existence of complete families with rho-value $1$ for at least one value of $D$.   
\end{remark}
\begin{appendix}

\begin{sidewaystable}
\vskip12cm   %1cm pour ce Mac, 12cm pour arXiv
\begin{center}
\caption{Values of $N(k,D,1.7,10^6,85698768)$ for $3\leq k\leq 30$ and various $D$ (see \S~\ref{numevid} for explanations)}\label{tablepreone}
\begin{tabular}{|c||c|c|c|c|c|c|c|c|c|c|c|c|c|c|c|}
\hline
   $D$& $1$&     $2$&     $3$&     $5$&    $6$&      $7$&    $10$&     $11$&   $13$&    $14$&   $15$&  $19$&     $23$&   $43$& $47$\\ 
\hline
 $h_D$& $1$&     $1$&     $1$&     $2$&    $2$&      $1$&    $2$&       $1$&    $2$&    $4$&     $2$&   $1$&     $3$&    $1$&   $5$\\
\hline
 $I_0$& $2000$& $1000$&  $3000$&  $500$&   $500$&    $1000$& $500$&    $1000$&   $500$& $250$&  $500$&  $1000$&  $333.3$& $1000$ &$200$\\
\hline
\hline
 $k=3$& $2087$& $1053$&   $0^*$&  $534$&   $512$&    $1012$& $514$&   $1049$&    $512$&  $246$&  $529$& $1049$& $362$&  $991$&  $195$\\
\hline
  $4$&   $0^*$&   $998$&  $3132$&  $568$&   $568$&   $1033$& $515$&    $1066$&   $510$&  $282$&  $507$& $1085$& $328$&  $992$&  $220$\\
\hline
  $5$&  $2193$&   $1001$&  $3219$&  $513$&  $544$&   $963$&  $552$&    $1079$&    $510$&  $271$& $507$& $1004$& $345$&  $1066$& $194$\\
\hline
  $6$&  $2118$&   $1008$&  $0^*$&  $535$&  $517$&   $1049$& $497$&    $1032$& $521$&   $261$&   $509$& $1088$&  $323$&  $1044$& $209$\\
\hline
  $7$&  $2107$&   $1024$&  $3112$&  $533$&   $517$& $2098^*$&$512$&    $1047$& $530$&  $270$&   $533$& $1061$&  $346$&  $1036$&  $208$\\
\hline
  $8$&  $4226^*$& $2117^*$& $3115$& $505$&   $520$&  $1018$&  $510$&   $1039$& $507$&   $249$&  $515$& $1056$&  $338$&  $1062$&  $174$\\
\hline
  $9$&  $2120$&   $1014$&  $6139^*$& $484$&   $503$&  $1041$&  $507$&    $984$& $512$&  $228$&  $549$& $1077$&  $329$&  $1060$&  $191$\\
\hline
 $10$&   $2167$&   $1039$&  $3171$&  $492$&   $536$&  $995$&   $509$&    $1038$& $539$&  $267$&  $523$& $990$&   $347$&  $1029$& $195$\\
\hline
$11$&   $2064$&   $1033$&  $3121$&   $518$&   $489$&  $1009$&  $447$&  $2084^*$& $524$&   $264$& $537$& $1035$&  $345$&  $1069$&  $205$\\
\hline
$12$&   $4239^*$&  $1048$&  $6368^*$& $519$& $547$&   $1009$&  $518$&   $1055$&  $502$&    $259$&  $519$& $1030$& $334$&  $1078$& $205$\\
\hline
$13$&   $1970$&    $1065$&  $3061$& $544$&  $504$&   $988$&   $476$&  $1059$&  $521$&    $229$&  $526$& $1076$& $333$&  $1028$& $192$\\
\hline
$14$&   $2095$&    $1102$&  $3243$& $560$&   $546$&   $2001^*$& $540$&  $1023$&  $532$&    $278$& $533$& $1048$& $364$&  $999$&  $225$\\
\hline
$15$&   $2030$&    $981$&   $6221^*$&  $526$&  $516$& $1130$&   $525$&    $982$& $502$&    $289$& $975^*$&$1058$& $347$&  $1077$&  $191$\\
\hline
$16$&   $4183^*$& $2058^*$& $3007$&  $528$& $536$& $1071$&  $502$&    $998$&    $511$&   $260$&  $491$& $1001$&  $361$&  $1071$& $205$\\
\hline
$17$&   $2073$&    $1008$&  $3194$&  $517$& $506$&  $1023$& $509$&    $1015$&   $482$&   $254$&  $470$& $1096$&   $374$& $1020$&  $206$\\
\hline
$18$&   $2139$&    $1017$& $6215^*$& $534$&  $512$&  $1013$& $537$&   $1021$&   $558$&   $273$&  $520$& $1016$&   $334$& $1001$&  $207$\\
\hline
$19$&   $2073$&    $1031$&  $3115$&  $529$&   $564$& $1049$& $497$&   $1048$&   $566$&   $229$&   $518$& $2127^*$& $356$&  $1025$& $205$\\
\hline
$20$&   $4063^*$&  $1071$&  $3111$&  $1073^*$& $517$& $1039$& $502$&  $1096$&   $481$&   $234$&   $491$& $1028$&   $325$&  $1101$& $196$\\
\hline
$21$&   $2035$&    $1068$&  $6304^*$& $526$&  $509$& $2016^*$& $500$&  $995$&   $568$&   $293$&   $503$& $1060$&   $371$&  $1019$&  $199$\\
\hline
$22$&   $2145$&    $996$&   $3048$&  $557$&   $512$& $1042$&  $533$&   $2138^*$& $519$&  $239$&   $545$& $1059$&   $345$&  $988$&  $216$\\
\hline
$23$&   $2113$&  $1012$& $3185$&   $530$&  $521$&   $1043$&  $476$&   $1071$&  $492$&   $271$&  $527$& $1059$& $682^*$& $1064$& $219$\\
\hline
$24$&  $4161^*$& $2110^*$&$6247^*$& $510$& $1055^*$& $1003$& $543$&   $996$&   $529$&   $260$&  $525$& $1031$& $333$&  $1113$&  $214$\\
\hline
$25$&   $1971$&  $1102$& $3082$&   $499$&  $504$&   $1031$&  $481$&   $1038$&  $540$&   $248$&   $523$&  $996$& $374$& $997$&  $227$\\
\hline
$26$&   $2065$&  $1055$& $3230$&  $493$&   $525$&   $1058$&  $542$&    $1042$&  $530$&   $257$&   $541$& $1083$& $336$& $1071$& $196$\\
\hline
$27$&   $2148$&  $1049$& $6327^*$&  $483$& $521$&   $1035$&  $516$&    $1062$&  $503$&   $270$&   $541$&  $976$& $323$& $1053$& $179$\\
\hline
$28$&   $4189^*$& $1038$& $3119$&  $547$&  $514$&   $2047^*$& $513$&  $1042$&  $506$&   $268$&  $480$&  $1006$& $367$&  $1054$&  $197$\\
\hline
$29$&   $2153$& $979$&  $3017$&   $581$&  $509$&   $1072$&   $551$&    $1040$& $522$&   $263$&  $500$&  $1030$& $334$&  $1086$&  $201$\\
\hline
$30$&   $2153$& $1041$& $6198^*$& $494$&  $535$&   $1029$&   $519$&    $1030$& $534$&   $271$&  $996^*$& $1068$& $361$& $955$&  $211$\\
\hline
\hline
Avg&    $2094.4$&$1034.8$& $3126.6$&  $524.8$& $522.6$& $1029.9$&$513.3$&  $1037.8$& $520.4$&$260.1$& $516.0$&      $1043.9$&   $345.6$&  $1041.0$&  $202.9$\\
\hline
\end{tabular}
\end{center}
\end{sidewaystable}

\begin{sidewaystable}
\vskip12cm  %1cm pour ce Mac, 12cm pour arXiv
\begin{center}
\caption{Values of $N(k,D,1.5,10^6,2\times 10^8)$ for $3\leq k\leq 30$ and various $D$ (see \S~\ref{numevid} for explanations)}\label{tableone}
\begin{tabular}{|c||c|c|c|c|c|c|c|c|c|c|c|c|c|c|c|}
\hline
   $D$& $1$&     $2$&     $3$&       $5$&    $6$&      $7$&    $10$&     $11$&   $13$&    $14$&   $15$&  $19$&   $23$&  $43$&  $47$\\ 
\hline
 $h_D$& $1$&     $1$&     $1$&      $2$&     $2$&      $1$&    $2$&       $1$&    $2$&    $4$&     $2$&   $1$&    $3$&   $1$&   $5$\\
\hline
 $I_0$& $116.3$& $58.17$& $174.5$&  $29.09$& $29.09$&  $58.17$&$29.09$&  $58.17$&$29.09$&$14.54$& $29.09$&$58.17$&$19.39$&$58.17$&$11.63$\\
\hline
\hline
  $k=3$& $132$&   $69$&    $0^*$&      $29$&   $34$&    $57$&    $35$&    $54$&    $29$&    $14$&   $27$&   $59$&   $17$&  $54$&  $12$\\
\hline
  $4$&   $0^*$&   $63$&    $198$&      $20$&   $31$&    $65$&   $31$&    $65$&     $27$&    $17$&   $37$&   $64$&   $22$&  $59$&  $10$\\
\hline
  $5$&   $123$&   $49$&    $211$&      $31$&   $26$&    $55$&   $24$&    $53$&    $30$&    $18$&   $26$&   $45$&   $21$&  $73$&  $12$\\
\hline
  $6$&   $132$&   $58$&    $0^*$&       $36$&  $41$&    $61$&    $22$&    $61$&    $32$&    $10$&   $29$&   $63$&   $14$&  $56$&  $13$\\
\hline
  $7$&   $111$&   $59$&    $190$&      $34$&   $32$&     $119^*$&  $29$&  $67$&    $32$&    $21$&   $27$&   $75$&   $15$&  $63$&  $6$\\
\hline
  $8$&   $235^*$&  $131^*$& $181$&   $30$&   $26$&     $56$&    $27$&    $47$&    $34$&    $16$&   $30$&   $64$&   $9$&  $61$&  $9$\\
\hline
  $9$&   $132$&   $60$&    $367^*$&   $31$&   $27$&     $52$&    $32$&    $63$&    $34$&    $22$&   $32$&   $80$&  $18$&  $52$&  $6$\\
\hline
 $10$&   $118$&   $55$&    $205$&     $28$&   $33$&     $69$&    $39$&    $59$&    $38$&    $13$&   $37$&   $46$&   $15$&  $66$&  $10$\\
\hline
$11$&   $111$&   $64$&    $197$&      $31$&   $38$&     $58$&    $26$&    $119^*$&    $29$&    $17$&   $28$&  $58$&   $15$&  $59$&  $13$\\
\hline
$12$&   $255^*$&  $42$&  $419^*$&    $22$&    $21$&     $62$&   $30$&    $67$&    $25$&    $27$&    $28$&      $61$&   $15$&  $59$&  $16$\\
\hline
$13$&   $125$&    $66$&    $164$&     $21$&    $27$&     $37$&   $26$&    $61$&    $43$&    $20$&    $32$&    $51$&   $28$&  $58$&  $9$\\
\hline
$14$&   $122$&    $74$&    $168$&     $29$&    $35$&     $133^*$&   $29$&  $45$&    $31$&    $13$&   $32$&   $55$&   $14$&  $69$&  $16$\\
\hline
$15$&   $119$&    $59$&    $381^*$&     $32$&    $30$&    $64$&   $28$&    $57$&    $30$&    $19$&    $57^*$&      $58$&   $16$&  $61$&  $9$\\
\hline
$16$&   $244^*$&    $130^*$&    $193$&     $30$&    $32$&     $58$&   $33$&    $53$&    $28$&   $9$&    $27$&      $71$&   $18$&  $77$&  $15$\\
\hline
$17$&   $133$&    $62$&    $194$&     $32$&    $33$&     $60$&   $22$&    $55$&    $30$&   $10$&    $36$&      $78$&   $16$&  $66$&  $11$\\
\hline
$18$&   $133$&    $59$&    $316^*$&     $34$&    $36$&     $65$&   $32$&    $62$&    $33$&   $18$&    $23$&    $63$&   $15$&  $71$&  $11$\\
\hline
$19$&   $111$&    $64$&    $176$&     $36$&    $27$&     $53$&   $31$&    $46$&    $38$&   $18$&    $32$&      $127^*$&   $24$&  $63$&  $15$\\
\hline
$20$&   $249^*$&    $60$&    $176$&     $64^*$&    $31$&     $73$&   $27$&    $57$&    $28$&   $12$&    $30$&      $63$&   $21$&  $61$&  $9$\\
\hline
$21$&   $113$&    $66$&    $378^*$&     $26$&    $25$&     $114^*$&   $26$&    $51$&  $33$&   $18$&    $30$&  $60$&   $25$&  $57$&  $12$\\
\hline
$22$&   $123$&    $62$&    $184$&     $25$&    $34$&     $55$&   $30$&    $127^*$&    $36$&   $19$&    $29$&  $68$&   $17$&  $54$&  $15$\\
\hline
$23$&   $103$&    $61$&    $192$&     $30$&    $44$&     $53$&   $38$&    $71$&    $32$&   $24$&    $17$&   $60$&   $44^*$&  $71$& $13$\\
\hline
$24$&   $207^*$& $129^*$&   $343^*$&   $28$&  $48^*$&    $64$&   $25$&    $69$&    $26$&   $14$&    $40$&   $60$&   $15$&  $51$&  $15$\\
\hline
$25$&   $96$&    $65$&    $186$&     $40$&    $26$&     $60$&   $33$&    $79$&    $34$&   $12$&    $28$&    $67$&   $20$&  $57$&  $10$\\
\hline
$26$&   $144$&    $57$&    $173$&   $33$&    $35$&     $66$&   $36$&    $65$&    $31$&   $14$&    $32$&      $45$&   $18$&  $59$&  $11$\\
\hline
$27$&   $135$&    $51$&    $354^*$&     $44$&    $40$&     $59$&   $27$&    $76$&    $21$&   $17$&    $17$&      $62$&   $27$&  $56$&  $10$\\
\hline
$28$&   $266^*$&    $66$&    $220$&     $25$&    $30$&   $123^*$&   $31$&    $66$&    $31$&   $19$&    $34$&    $65$&   $23$&  $71$&  $11$\\
\hline
$29$&   $113$&    $69$&    $170$&     $34$&    $23$&     $69$&   $29$&    $60$&    $26$&   $21$&    $25$&      $69$&   $23$&  $43$&  $12$\\
\hline
$30$&   $109$&    $67$&    $388^*$&     $24$&    $37$&     $47$&   $26$&    $55$&    $29$&   $13$&    $69^*$&      $47$&   $25$&  $55$&  $12$\\
\hline
\hline
Avg&    $121.0$&$61.50$& $186.6$&  $30.25$& $31.36$& $59.38$&$29.43$&  $60.25$& $31.07$&$16.61$& $29.57$&      $61.45$&   $18.86$&  $60.79$&  $11.54$\\
\hline
\end{tabular}
\end{center}
\end{sidewaystable}

\newpage

\begin{center}  
Data for $k=28$, $D=1$, $\rho_0\in \{1.1, 1.2,1.3,1.4,1.5\}$.
\vskip1mm
\begin{tabular}{|c||c|c|c|}
\hline
Interval& $10^6\leq r\leq 10^8$& $10^8\leq r\leq 10^{10}$& $10^{12}-10^{10}\leq r\leq 10^{12}+10^{10}$\\
\hline
\hline
$N(1.1)$& $0$                  &    $1$                  &   $0$  \\
\hline
$I(1.1)$& $0.325$              &   $0.311$               &  $0.002$ \\
\hline
$N(1.2)$& $3$                  &    $6$                  &   $0$   \\
\hline
$I(1.2)$& $1.502$              &   $2.286$               &   $0.022$ \\
\hline
$N(1.3)$& $8$                  &    $24$                 &   $0$   \\
\hline
$I(1.3)$& $7.104$              &   $17.22$               &   $0.321$ \\
\hline
$N(1.4)$& $37$                 &    $135$                &   $5$   \\
\hline
$I(1.4)$& $34.39$              &   $132.71$              &   $4.723$  \\
\hline
$N(1.5)$& $188$                &    $1128$               &   $73$   \\
\hline
$I(1.5)$& $170.07$             &    $1044.7$             &   $69.86$ \\
\hline                
\end{tabular}
\end{center}
\medskip

\begin{center}  
Data for $k=27$, $D=11$, $\rho_0\in \{1.1, 1.2,1.3,1.4,1.5\}$.
\vskip1mm
\begin{tabular}{|c||c|c|c|}
\hline
Interval& $10^6\leq r\leq 10^8$& $10^8\leq r\leq 10^{10}$& $10^{12}-10^{10}\leq r\leq 10^{12}+10^{10}$\\
\hline
\hline
$N(1.1)$& $0$                  &    $0$                  &   $0$  \\
\hline
$I(1.1)$& $0.081$              &   $0.078$               &  $0.00038$ \\
\hline
$N(1.2)$& $0$                  &    $2$                  &   $0$   \\
\hline
$I(1.2)$& $0.375$              &   $0.57$               &   $0.0055$ \\
\hline
$N(1.3)$& $1$                  &    $5$                 &   $0$   \\
\hline
$I(1.3)$& $1.78$              &   $4.31$               &   $0.080$ \\
\hline
$N(1.4)$& $9$                 &    $30$                &   $1$   \\
\hline
$I(1.4)$& $8.60$              &   $33.18$              &   $1.18$  \\
\hline
$N(1.5)$& $57$                &    $271$               &   $22$   \\
\hline
$I(1.5)$& $42.52$             &    $261.17$             &   $17.46$ \\
\hline                
\end{tabular}
\end{center}
\medskip

\begin{center}  
Data for $k=8$, $D=23$, $\rho_0\in \{1.1, 1.2,1.3,1.4,1.5\}$.
\vskip1mm
\begin{tabular}{|c||c|c|c|}
\hline
Interval& $10^6\leq r\leq 10^8$& $10^8\leq r\leq 10^{10}$& $10^{12}-10^{10}\leq r\leq 10^{12}+10^{10}$\\
\hline
\hline
$N(1.1)$& $0$                  &    $0$                  &   $0$  \\
\hline
$I(1.1)$& $0.027$              &   $0.026$               &  $0.00013$ \\
\hline
$N(1.2)$& $0$                  &    $0$                  &   $0$   \\
\hline
$I(1.2)$& $0.125$              &   $0.191$               &   $0.00183$ \\
\hline
$N(1.3)$& $0$                  &    $1$                 &   $0$   \\
\hline
$I(1.3)$& $0.592$              &   $1.435$               &   $0.0267$ \\
\hline
$N(1.4)$& $1$                 &    $16$                &   $0$   \\
\hline
$I(1.4)$& $2.866$              &   $11.06$              &   $0.394$  \\
\hline
$N(1.5)$& $7$                &    $76$               &   $6$   \\
\hline
$I(1.5)$& $14.17$             &    $87.06$             &   $5.821$ \\
\hline                
\end{tabular}
\end{center}

\end{appendix}

\end{document}